\documentclass[a4paper,11pt]{amsart}
\usepackage{amssymb}
\usepackage{amsmath}
\usepackage{mathrsfs}
\usepackage{amsmath,amssymb,amsthm,latexsym,amscd,mathrsfs}
\usepackage{indentfirst}
\usepackage{graphicx}
\usepackage{booktabs}
\usepackage{array}
\usepackage{chemarrow}
\newcommand{\ROM}[1]{\mathrm{\uppercase\expandafter{\romannumeral#1}}}
\theoremstyle{definition}

\newtheorem{thm}{Theorem}[section]

\newtheorem{cor}{Corollary}[section]

\newtheorem{rem}{Remark}[section]

\newtheorem{ack}{Acknowledgements}   
\makeatletter


\title{$\eta$-invariant and a problem of B\'erard-Bergery on the existence of closed geodesics}
\author[Zizhou Tang]{Zizhou Tang}\address{School of Mathematical Sciences, Laboratory of Mathematics and Complex Systems, Beijing Normal
University, Beijing 100875, P.R.China}\email{zztang@bnu.edu.cn}
\thanks {The project is partially supported by MOEC and NNSFC.}
\author[Weiping Zhang]{Weiping Zhang}\address{Chern Institute of Mathematics and LPMC,
Nankai University, Tianjin 300071,
P.R.China}\email{weiping@nankai.edu.cn}

\date{}
\begin{document}

\maketitle
\begin{abstract} We use the $\eta$-invariant of Atiyah-Patodi-Singer to compute the Eells-Kuiper invariant for the 
  Eells-Kuiper quaternionic projective plane. By combining with a known result of  B\'{e}rard-Bergery, it shows that every  Eells-Kuiper quaternionic projective plane   carries a Riemannian metric such that all geodesics passing through a certain point are simply closed and of the same length. 
\end{abstract}


\section{Introduction}

The $\eta$-invariant introduced by Atiyah-Patodi-Singer \cite{APS}, as well as its various ramifications,  has played important roles in many prblems in geometry and topology.  In this short paper, we   use the $\eta$-invariant to compute the Eells-Kuiper invariant for the 
  Eells-Kuiper quaternionic projective plane. By combining with a known result of  B\'{e}rard-Bergery, it shows that every  Eells-Kuiper quaternionic projective plane   carries a Riemannian metric such that all geodesics passing through a certain point are simply closed and of the same length. 

  To be more precise, let $p$ be a point  in a closed manifold $M$. Let $g$ be a Riemannian metric on $M$. The Riemannian structure $(M,g)$ is called an $SC^p$
Riemannian structure if all  geodesics issued from $p$ are
simply  closed (periodic) geodesics with the same length.  
We refer  to the classic  book \cite{Be} for a systematic acount of the $SC^p$ structures. 

  It is clear   that there are  $SC^p$ Riemannian structures on the compact symmetric spaces of rank one (briefed  in \cite{Be} as CROSS), namely  the
unit spheres, the real projective spaces, the complex projective
spaces, the quaternionic projective spaces and the Cayley projective
plane, endowed with the corresponding  canonical metrics. Moreover, a fundamental result of Bott \cite{Bo}  states  that any smooth manifold carrying  an $SC^p$
structure should have the same integral cohomolgy ring as that of a 
CROSS. On the other hand, there are manifolds verifying the above cohomological condition but not diffeomorphic to any CROSS.  For typical examples, we mention   the (exotic) homotopy spheres and the Eells-Kuiper (exotic) quaternionic
projective planes. 


In 1975, B\'{e}rard-Bergery \cite{BB} discovered  an $SC^p$ structure on an exotic
sphere of dimension $10$. He then raised   the natural question:
\emph{is there any (exotic) Eells-Kuiper quaternionic projective plane carrying 
an $SC^p$ structure}?  The same question was also posed explicitly by Besse in the classic book    \cite[0.15 on pp. 4]{Be}.
 Moreover, it is pointed out   in \cite[pp. 143]{Be}  that a positive answer to the above question would also give 
a positive nontrivial example  to the following open question: 
   {\it whether  
a   Blaschke manifold at a point\footnote{Cf. \cite[5.37 on page 135]{Be} for a definition.} would carry an  $SC^p$ Riemannian structure}? 

The purpose of this article is to provide a positive answer to the above two questions concerning the Eells-Kuiper  quaternionic projective planes. 

Before going on,  we describe  the Eells-Kuiper  quaternionic projective planes as follows,
starting  with  the  standard  construction of 
 Milnor \cite{Mi1}. 

For any  pair of integers $(h,j)$, 
let  $\xi_{h,j}$ be the $S^3$-bundle over $S^4$ 
  determined by the
characteristic map $f_{h,j}: S^3 \longrightarrow SO(4)$ with $
f_{h,j}(u)v=u^hvu^j$ for $u \in S^3$, $v \in {\bf R}^4$,  where we identify
${\bf R}^4$ with the space of quaternions. 
 It is shown in \cite{Mi1} that   when $h+j=1$, the total space of the above 
sphere bundle  is homeomorphic to the unit sphere $S^7$.  From now on, we denote by $M_h$  this 
total space corresponding to $(h,j)=(h,1-h)$, and denote by $N_h$  the
associated disk  bundle.

\begin{rem}
When $h=0$ or $1$, $M_h$ is just the unit 7-sphere and the sphere
bundle is just the Hopf fibration (corresponding to the left or
right multiplications of the quaternions, respectively). On the other hand,
$M_2$ is the exotic sphere generating the group $\Theta(7)$ (the set of
the orientation preserving diffeomorphism classes of $7$-dimensional
oriented homotopy spheres), which is isomorphic to the cyclic group
${\bf Z}_{28}$.
\end{rem}

It is shown by Eells-Kuiper \cite{EK2} that the homotopy sphere $M_h$ is
diffeomorphic to $S^7$ if and only if the following congruence holds for $h$, 
             $$\frac{ h(h-1)}{56}\equiv 0 \ \  {\rm mod}\ \ {\bf Z}. \eqno(1.1)$$

From now on, we assume that $h$ satisfies (1.1).
Then there is a diffeomorphism  $\sigma:M_h\rightarrow S^7$. Let $X_{h,\sigma}$ denote the $8$ dimensional closed smooth
manifold constructed  from $N_h$ by attaching the unit disk  $D^8$ by the
diffeomorphism $\sigma:\partial(N_h)=M_h\rightarrow \partial(D^8)=S^7$. This is what we call an
Eells-Kuiper  quaternionic projective plane, first constructed in \cite{EK1}.\footnote{Indeed, Eells and Kuiper showed in \cite{EK1} that the $X_{h,\sigma}$'s are the only $8$ dimensional closed smooth manifolds  admitting a Morse function with $3$ critical points. } We
remark   that when $h=0$ or $1$, and $\sigma={\rm id}$,  
$X_{h,\sigma}$ is just the standard quaternionic projective plane
${\bf H}P^2$.
We also mention a  deep result due to Kramer and   Stolz \cite{KS1} which 
asserts that the diffeomorphism type of the resulting manifold
$X_{h,\sigma}$ does not depend on the choice of the diffeomorphism
$\sigma:M_h\rightarrow S^7$.

Let $\tau_h$ be the canonical involution on $M_h$ obtained by the fiberwise antipodal involution on $S^3$. 
By \cite[Theorem 1]{BB} and the above result of Kramer-Stolz, to prove that $X_{h,\sigma}$ carries an $SC^p$ Riemannian structure, one only needs to show that
there is a diffeomorphism $\sigma':M_h\rightarrow S^7$ such that $ \tau \sigma'=\sigma' \tau_h$, where $\tau$ is the standard antipodal involution of $S^7$. Equivalently, one needs only to show that 
 the quotient manifold $M_h/\tau_h$ is diffeomorphic to 
${\bf R}P^7$.  This is the content of the following main result of this paper. 

\begin{thm}\label{thm}
\emph{ The involution $\tau_h$ on $M_h\cong S^7$ is equivalent to
the standard antipodal involution on $S^7$. In other words,  
$M_h/\tau_h$ is diffeomorphic to ${\bf R}P^7$. }
\end{thm}

\begin{cor}\label{cor}
\emph{ Every  Eells-Kuiper quaternionic projective plane
admits an $SC^p$ Riemannian structure.}
\end{cor}

 \begin{rem}
Since there is infinitely many Eells-Kuiper quaternionic projective planes not diffeomorphic to each other, the above Corollary actually shows that there is an infinite family of pairwise non-diffeomorphic  manifolds $M$ with the cohomology ring of ${\bf H}P^2$ such that each $M$ admits    an $SC^p$ Riemannian structure. 
\end{rem}

The rest of this article is organized as follows. In Section 2, we reduce the proof of Theorem 1.1 to a problem of computing the Eells-Kuiper $\mu$ invariant introduced in \cite{EK2}. In Section 3, we recall the results of Donnelly \cite{Do1} and Kreck-Stolz \cite{KS2} (cf. \cite{G}) which use $\eta$-invariants to express the $\mu$-nvariant, and then  carry out the required computation of the involved $\eta $ invariant.

\section{Theorem 1.1 and the  Eells-Kuiper $\mu$ invariant}

As was indicated in \cite[pp. 240]{BB}, by results of Mayer \cite{Ma}, there could  only be two possibilities for  $M_h/\tau_h$. That is, it is diffeomorphic either  to ${\bf R}P^7$ or to the connected sum  ${\bf R}P^7\# 14 M_2$, where $14M_2$ is the connected sum $M_2\#\cdots\# M_2$ of $14$ copies of $M_2$. 

On the other hand, Milnor \cite{Mi2} showed that  the  Eells-Kuiper $\mu$ invariant of ${\bf R}P^7$ and  ${\bf R}P^7\# 14 M_2$ takes different values. Thus, in order to prove Theorem 1.1, one need only to show that the $\mu$ invariant of $M_h/\tau_h$ is different from that of  ${\bf R}P^7\# 14 M_2$.

For completeness, we recall the definition of the Eells-Kuiper $\mu$ invariant in our situation. 
Let $M$ be a $7$ dimensional closed oriented spin manifold such that the $4$-th cohomology group $H^4(M;{\bf R})$ vanishes.\footnote{By the above diffeomorphism type result, it is clear that $M_h/\tau_h$ verifies this condition.}  If $M$ bounds a compact  oriented spin manifold $N$, then the first Pontrjagin class $p_1(N)\in H^4(N,M;{\bf Q})$ is well-defined.

Following 
\cite[(11)]{EK2}, we  define $\mu(M)\in{\bf R}/{\bf Z}$ by 
$$ \mu(M)\equiv\frac{p^2_1(N)}{2^7\times 7}-\frac{{\rm Sign}(N)}{2^5\times 7}\ \ {\rm mod}\ \ {\bf Z},\eqno(2.1)$$
where $p_1^2(N)$ denotes the corresponding Pontrjagin number and  ${\rm Sign}(N)$ is the
Signature of  $N$. 

Now set $M=M_h$, $N=N_h$. 
Let  $x\in H^4(S^4;{\bf Z})$ be 
the generator.
By \cite{Mi1}, one has
$$ e(\xi_{h,1-h})=x,\ \ \ \  p_1(\xi_{h,1-h})=\pm 2(2h-1)x,\eqno(2.2) $$
where $e(\xi_{h,1-h})$ and $p_1(\xi_{h,1-h})$ are the Euler class and
the first Pontrjagin class of the sphere bundle $\xi_{h,1-h}$ respectively. Also by \cite{Mi1}, one has
$${\rm Sign}(N_h)=1.\eqno(2.3)$$

From (2.2) and (2.3), one deduces as in \cite{Mi1} and \cite{EK2} that
$$\frac{p^2_1(N_h)}{2^7\times 7}-\frac{{\rm Sign}(N_h)}{2^5\times 7}=\frac{h(h-1)}{56},\eqno(2.4)$$
which is an integer in view of the assumption (1.1). 

Recall that  by \cite{Mi2}, one has  $\mu({\bf R}P^7)=\pm \frac{1}{32}$ while $\mu( {\bf R}P^7\# 14 M_2)=\pm\frac{1}{32}+\frac{1}{2}$. Thus, in order to prove Theorem 1.1, one need only to prove the following result. 

\begin{thm}\label{thm}  {\it The following identity holds for any integer $h$ verifying (1.1),}
 $$\mu(M_h/\tau_h)\equiv \pm \frac{1}{32}\ \ \ {\rm mod}\ \ {\bf Z}.\eqno(2.5)$$
\end{thm}

Theorem 2.1 will be proved in Section 3

\section{A proof of Theorem 2.1}

In this section,   we compute   $\mu(M_h/\tau_h)$. The obvious  difficulty   is that one does not   find easily an $8$ dimensional spin manifold   with boundary $M_h/\tau_h$. Instead, we will make use of an intrinsic formula for the $\mu$ invariant, which is given by Donnelly \cite{Do1} and Kreck-Stolz \cite{KS2} (cf. the survey paper of Goette \cite{G}). 

Indeed, for any $7$ dimensional closed oriented spin manifold $M$ with $H^4(M;{\bf R})=0$, let $g^{TM}$ be a Riemannian metric on $TM$. Let $\nabla^{TM}$ be the associated Levi-Civita connection. Let $p_1(TM,\nabla^{TM})$ be the corresponding first Pontrjagin form (cf.  \cite[Section 1.6.2]{Z}). Then there is a   $3$-form $\widehat{p}_1(TM,\nabla^{TM})$ on $M$ such that
$$d \,
\widehat{p}_1(TM,\nabla^{TM})= p_1(TM,\nabla^{TM}).\eqno(3.1)$$

Let $D_M$ (resp. $B_M$) be the  Dirac (resp.
Signature)  operator associated to $g^{TM}$. Let $\eta(D_M)$, $\eta(B_M)$ be the Atiyah-Patodi-Singer $\eta$ invariant of $D_M$, $B_M$   (cf. \cite{APS}).  Let  $$\overline{\eta}(D_M)=\frac{1}{2}\left(\dim(\ker D_M)
+\eta(D_M)\right)$$ be the corresponding  reduced $\eta$-invariant. 
 
By \cite{Do1} and  \cite{KS2} (cf. \cite[pp. 424]{G}),  the  $\mu$ 
invariant defined in (2.1) can be represented  by
$$\mu(M)\equiv \overline{\eta}(D_M)+\frac{\eta(B_M)}{2^5\times 7}-\frac{1}{2^7\times 7}\int_{M} p_1(TM,\nabla^{TM})\wedge \widehat{p}_1(TM,\nabla^{TM})\ \ {\rm mod}\ \ {\bf Z}.\eqno (3.2)$$

Now consider the double covering $M_h\rightarrow M_h/\tau_h$. We fix a spin structure on $M_h/\tau_h$ and
lift everything from $M_h/\tau_h$ to $M_h$. We get that
$$\mu(M_h /\tau_h )\equiv \overline{\eta}(P_hD_{M_h})+\frac{\eta(P_hB_{M_h})}{2^5\times 7}-\frac{1}{2^8\times 7}\int_{M_h} p_1(TM_h,\nabla^{TM_h})\wedge \widehat{p}_1(TM_h,\nabla^{TM_h})\  {\rm mod}\  {\bf Z},\eqno (3.3)$$
where $P_h=\frac{1}{2}(1+\tau_h)$ is the canonical projection. Here $\tau_h$ denotes the lifted actions on the corresponding vector bundles.

Indeed, recall that $M_h$ is a fiber bundle over $S^4$ with  fiber $S^3$. It is the boundary of the unit  disk bundle $N_h$ over $S^4$, while $\tau_h$ is the canonical involution which maps on each fiber by mapping a point to its antipodal. This involution  extends canonically to an   involution on $N_h$ which we still denote by $\tau_h$.
Clearly, the fixed point set of $\tau_h$ on $N_h$ is $S^4$, the image of the zero section of the disk bundle.

Let $g^{TN_h}$ be a $\tau_h$ invariant Riemannian metric on $TN_h$ such that it restricts to $g^{TM_h}$ on $\partial N_h=M_h$ and   is of product structure near $ M_h$ (the existence of such a metric is clear).  Let $\nabla^{TN_h}$ be the associated Levi-Civita connection. 

By dimensional reason we see that we are in  the situation of even type in the sense of   \cite[Proposition 8.46]{AB}.  Thus
there exists a  $\tau_h$-equivariant spin structure on $N_h$, such that it induces a $\tau_h$-equivariant spin structure on $M_h$, which  equals  the one  lifted from the spin structure given on $M_h/\tau_h$.   In particular, $\tau_h$ lifts to an action on the associated spinor bundle $S(TN_h)=S_+(TN_h)\oplus S_-(TN_h)$ associated to $(TN_h,g^{TN_h})$, preserving the corresponding ${\bf Z}_2$-grading. It induces an action on $S(TM_h)=S_+(TN_h)|_{M_h}$. Moreover, the lifted $\tau_h$ action  commutes with the   Dirac operator $D_{N_h}:\Gamma(S(TN_h))\rightarrow \Gamma(S(TN_h))$, and thus also commutes with the induced Dirac operator $D_{M_h}:\Gamma(S(TM_h))\rightarrow \Gamma(S(TM_h))$, which in turn determines a Dirac operator on $M_h/\tau_h$ on which one can apply (3.2) and (3.3). 

 Let  $D_{N_h,+}:\Gamma(S_+(TN_h))\rightarrow \Gamma(S_-(TN_h))$ be the natural restriction   of $D_{N_h}$. 
By  the Atiyah-Patodi-Singer index theorem \cite{APS} and its
equivariant extention by Donnelly \cite{Do2}, one finds
$$ \overline{\eta}(P_hD_{M_h}) \equiv \frac{1}{2}\int_{N_h}\widehat{A}(TN_h,\nabla^{TN_h}) + \frac{1}{2}\int_{S^4} A_1\ \ {\rm mod}\ \ {\bf Z},\eqno(3.4)$$
where the mod ${\bf Z}$ term comes from the Atiyah-Patodi-Singer type index ${\rm ind}_{\rm APS}(P_hD_{N_h,+})$, $\widehat{A}(TN_h,\nabla^{TN_h})$ is the Hirzebruch $\widehat{A}$-form associated to $\nabla^{TN_h}$ (cf. \cite[Section 1.6.3]{Z})  and  $A_1$ is the canonical contribution on the fixed point set (which by the local index theory is the same as the usual fixed point set contribution appearing in the equivariant Atiyah-Singer index theorem for compact group actions on closed manifolds). 

Similarly,
$$  {\eta}(P_hB_{M_h}) = \frac{1}{2}\int_{N_h}  {L}(TN_h, \nabla^{TN_h}) + \frac{1}{2}\int_{S^4} A_2 - \frac{1}{2} {\rm Sign}(N_h)-\frac{1}{2} {\rm Sign}(N_h,\tau_h) ,\eqno(3.5)$$
where  $ {L}(TN_h, \nabla^{TN_h}) $ is the Hirzebruch $L$-form associated to $\nabla^{TN_h}$ (cf. \cite[Section 1.6.3]{Z}), $A_2$ is the canonical contribution on the fixed point set and   ${\rm Sign}(N_h,\tau_h)$ is the notation for the equivariant
Signature with respect to $\tau_h$.  

By a direct computation, one has
$$ \frac{1}{2}\int_{N_h}\widehat{A}(TN_h,\nabla^{TN_h}) +\frac{1}{2^6\times 7} \left({\int_{N_h}  {L}(TN_h, \nabla^{TN_h}) - {\rm Sign}(N_h)} \right)$$
$$ -\frac{1}{2^8\times 7}\int_{M_h} p_1(TM_h,\nabla^{TM_h})\wedge \widehat{p}_1(TM_h,\nabla^{TM_h})=\frac{p^2_1(N_h)}{2^8\times 7}-\frac{{\rm Sign}(N_h)}{2^6\times 7}.\eqno(3.6)$$

From (2.4) and (3.3)-(3.6),  
we find that
$$\mu(M_h /\tau_h)\equiv \frac{h(h-1)}{112}  + \frac{1}{2}\int_{S^4} A_1+\frac{1}{2^6\times 7}\int_{S^4}A_2-\frac{{ \rm Sign}(N_h,\tau_h)}{2^6\times 7}\ \ {\rm mod}\ \ {\bf Z}.\eqno(3.7)$$

Now let $W_h$ denote the normal  bundle  in $N_h$ to the submanifold $S^4$,   the fixed point set of $\tau_h$.  It is clear that $\tau_h$ acts on   $W_h$ by multiplication by $-1$.

By (2.2) and \cite[pp. 267]{LM}, one finds 
$$\int_{S^4}A_1=\pm\frac{1}{32}\int_{S^4}p_1(W_h)=\pm \frac{(2h-1)}{16}.\eqno(3.8)$$
Similarly, by \cite[pp. 265]{LM} and (2.2), one has
$$\int_{S^4}A_2=\int_{S^4}e(W_h)=1.\eqno(3.9)$$

On the other hand, since $S^4$ is the fixed point set of $\tau_h$, $\tau_h$ preserves $x\in H^4(S^4;{\bf Z})$. Thus one has
$${\rm Sign}(N_h,\tau_h) =1.\eqno(3.10)$$

From (3.7)-(3.10), one gets
$$\mu(M_h /\tau_h)\equiv \frac{h(h-1)}{112}  \pm \frac{2h-1}{32}  \ \ {\rm mod}\ \ {\bf Z}.\eqno(3.11)$$

We now claim that under the condition (1.1), (2.5) follows from (3.11). 

Indeed, under the assumption (1.1), one has $ h\equiv 0,\, 1,\, 8,\, 49\ {\rm
mod}\  56{\bf Z}$. Thus we only need to do the case by case checking as follows, where by ``$\equiv$'' we mean that the congruence is  $   {\rm mod}\ {\bf Z}$.

(i) For $h=56k$, then $\frac{h(h-1)}{112}\equiv \frac{k}{2}$, while
$\frac{2h-1}{32}\equiv -\frac{1}{32}+\frac{k}{2}$;

(ii) For $h=56k+1$, then $\frac{h(h-1)}{112}\equiv \frac{k}{2}$,
while $\frac{2h-1}{32}\equiv \frac{1}{32}+\frac{k}{2}$;

(iii) For $h=56k+8$, one has $\frac{h(h-1)}{112}\equiv \frac{1}{2}
+\frac{k}{2}$, while $\frac{2h-1}{32}\equiv
-\frac{1}{32}+\frac{1}{2}+\frac{k}{2}$;

(iv) For $h=56k+49$, one has $\frac{h(h-1)}{112}\equiv \frac{k}{2}$,
while $\frac{2h-1}{32}\equiv \frac{1}{32} +\frac{k}{2}$.

Combining (i)-(iv) with (3.11), we always have (2.5).

The proof of Theorem 2.1, as well as of   Theorem 1.1 and Corollary 1.1 is complete.

\begin{ack}
The first author thanks Dr. Chao Qian for useful discussion. We also thank the referee of this paper  for very helpful suggestions. 
\end{ack}

\end{document}